\documentclass[a4paper, reqno,10pt]{amsart}
\usepackage{amsmath}
\usepackage{graphicx}
\usepackage{amsfonts}
\usepackage{amssymb}
\newtheorem{theorem}{Theorem}[section]
\newtheorem{lemma}[theorem]{Lemma}
\newtheorem{proposition}[theorem]{Proposition}

\newtheorem{hypothesis}[theorem]{Hypothesis}
\theoremstyle{definition}
\newtheorem{definition}[theorem]{Definition}

\theoremstyle{remark}

\numberwithin{equation}{section}
\def\Dim{\noindent\hbox{{\bf Proof.}$\;\; $}}
\def\finedim{{\hfill\hbox{\enspace${ \square}$}} \smallskip}
\def\sqr#1#2{{\vcenter{\vbox{\hrule height .#2pt
\hbox{\vrule width .#2pt height#1pt \kern#1pt \vrule
width .#2pt} \hrule height .#2pt}}}}
\def\square{\mathchoice\sqr54\sqr54\sqr{4.1}3\sqr{3.5}3}

\def\nat{\mathbb N}

\def\P{{\mathcal P}}
\def\Pr{\mathbb P}
\def\E{\mathbb E}

\def\R{\reali}

\def\reali{\mathbb R}

\def\F{{\mathcal F}}

\def\G{{\mathcal G}}
\def\P{{\mathcal P}}
\def\U{{\mathcal U}}

\def\Et{\mathbb{E}^{\F_t}}

\begin{document}
\title{Ergodic Optimal Quadratic Control for an Affine Equation with Stochastic and Stationary Coefficients}

\author{Giuseppina Guatteri and Federica Masiero}

\address{Giuseppina Guatteri, Dipartimento di Matematica,
  Politecnico di Milano, piazza Leonardo da Vinci 32, 20133 Milano,
e-mail: giuseppina.guatteri@polimi.it}
\address{Federica Masiero,
Dipartimento di Matematica e Applicazioni, Universit\`{a} di Milano
Bicocca, via R. Cozzi 53 - Edificio U5, 20125 Milano, e-mail:
federica.masiero@unimib.it}
\maketitle

\begin{abstract}
We study ergodic quadratic optimal stochastic control problems for
an affine state equation with state and control dependent noise
and with stochastic coefficients. We assume stationarity of the
coefficients and a finite cost condition.
We first treat the stationary case
and we show that the optimal cost corresponding to this ergodic control
problem coincides with the one corresponding to a suitable
stationary control problem and we provide a full
characterization of the ergodic optimal cost and control.
\end{abstract}



\medskip


\bibliographystyle{plain}

%
%
%
%
%
%
%
%
%
%
\section{Introduction}

In this paper we study an ergodic quadratic control problem for a linear affine
equation with both state and  control dependent noise, and the coefficients of the
state equation, allowed to be random, are assumed to be stationary. We continue
our previous work \cite{GM},
where the infinite horizon case and the ergodic case are studied but no
characterization of the ergodic limit was given. The main result of
the present paper is to obtain the characterization of the ergodic limit,
see Theorem \ref{teofinale}, when the coefficients are stationary in a
suitable sense, see \cite{Tess2} and section 2 below.

The main tool will be Backward Stochastic Riccati Equations (BSREs): such equations
are naturally linked with stochastic
optimal control problems with stochastic coefficients.
The first existence and uniqueness result for such a kind of equations
has been given by Bismut in \cite{Bi76}, but then several works, see e. g.
\cite{Bi78}, \cite{KoTa01}, \cite{KoTa02}, \cite{KoTa03}, \cite{Peng} followed.
Only very recently Tang in \cite{Tang} solved the general non singular
case corresponding to the
linear quadratic problem with random coefficients and control
dependent noise. In \cite{GM}, we have studied the infinite horizon case
and the ergodic case
namely, we have considered a cost functional depending only on the
asymptotic behaviour of the state (ergodic control).

Starting from this point, in this paper  we
first consider the stationary problem: minimize over
all admissible controls the cost functional
\begin{equation*}
J^{\natural}(u,X)=\E\int_{0}^{1}[\vert \sqrt{S_s}X_s \vert^2+\vert u_s\vert^2] ds.
\end{equation*}
The control $u$ is stationary and $X$ is the corresponding solution
of the state equation
\begin{equation}
dX_t=A_tX_tdt+B_tu_tdt+
\sum_{i=1}^dC^i_tX_tdW^i_t+\sum_{i=1}^dD^i_tu_tdW^i_t+f_tdt.
\end{equation}
We denote the optimal cost for the stationary problem by $\overline{J}^{\natural}$.

\noindent The main technical point of this paper is to prove that the
closed loop equation for the stationary control problem,
admits a unique stationary solution, see proposition
\ref{propclosedloopSTAZ}.

In order to study the ergodic control problem, we first consider the discounted cost
functional
\begin{equation}
J^{\alpha}(0,x,u)=\mathbb{E}\int_{0}^{+\infty}e^{-2\alpha s}[\left\langle S_{s}X_{s}%
^{0,x,u},X_{s}^{0,x,u}\right\rangle +|u_{s}|^{2}]ds,
\label{costo__intro}
\end{equation}
where $X$ is solution to equation
\begin{equation}
\left\{
\begin{array}
[c]{ll}
dX_{s}=(A_{s}X_{s}+B_{s}u_{s})ds+
{\displaystyle\sum_{i=1}^{d}}
\left(  C_{s}^{i}X_{s}+D_{s}^{i}u_{s}\right)  dW_{s}^{i} +f_s ds, & s\geq 0\\
X_{0}=x. &
\end{array}
\right.  \label{f.stato.stato.intro}%
\end{equation}
$A$, $B$, $C$ and $D$ are bounded random and stationary processes and $f\in
L^{\infty}_\P(\Omega\times[0,+\infty),\R^n)$, moreover we assume suitable finite cost
conditions. It is proved in \cite{GM} that in general, without stationarity
assumptions,
\begin{multline*}
 \underline{\lim}_{\alpha\rightarrow 0}\alpha \overline{J}^{\alpha}(x)=  \underline{\lim}_{\alpha\rightarrow 0}\alpha \E\int_0^{+\infty}2\langle r^{\alpha}_s, f^{\alpha}_s\rangle ds \\ -\overline{\lim}_{\alpha\rightarrow 0}\alpha\E\int_0 ^{+\infty} \vert (  I+
{\displaystyle\sum_{i=1}^{d}}
\left(  D_{s}^{i}\right)  ^{\ast}P^{\alpha}_{s}D_{s}^{i})  ^{-1}( B_{s}^{*}r^{\alpha}_{s}+{\displaystyle\sum_{i=1}^{d}}\left( D_{s}^{i}\right) ^{*}g_{s}^{\alpha,i})\vert ^2 ds.
\end{multline*}
Starting from this point, we show here that in the stationary case
\begin{equation*}
 \underline{\lim}_{\alpha\rightarrow 0}2\alpha \overline{J}^{\alpha}(x)= \overline{J}^{\natural}
\end{equation*}
Then we consider the ``true'' ergodic optimal cost, we minimize the following functional
\[
\widehat{J}(x,u)=\underline{\lim}_{\alpha\rightarrow 0} 2\alpha J(x,u)
\]
over all $u\in \widehat{\U}$, see (\ref{Ucappuccio}) for the
definition of $\widehat{\U}$.
We are able to prove that
\[
\inf_{u\in \widehat{\U}}\widehat{J}(x,u)=J^{\natural}(u).
\]
and to the characterize the optimal ergodic control, see lemma
\ref{lemmafinale} and theorem \ref{teofinale}.

\section{Linear Quadratic optimal control in the stationary case}

Let $(\Omega, \F,  \Pr)$ be a probability space and assume that $W:(-\infty, +\infty)\rightarrow \reali$ is a $d$-dimensional brownian motion defined on the whole real axis. Let $\{\F_t\}_{ t \in (-\infty,+\infty)}$ its natural filtration completed. For all $s,t \in \reali$ with $t\geq s$ we denote by $\G^s_t$ the $\sigma$-field generated by $\left\lbrace   W_{\tau}-W_s, s\leq \tau \leq t \right\rbrace$. Notice that for all $s \in \reali$, $\left\lbrace \G^s_t\right\rbrace_{t\geq s}$ is a filtration in $(\Omega, \F)$. Finally we assume that for all $s<0$, $\G^s_0\subseteq \F_0$.

\noindent Next we set a stationary framework:
we introduce the semigroup $(\theta_t)_{t\in \reali}$  of measurable
mappings  $\theta_t: (\Omega,\mathcal{E})\rightarrow
(\Omega,\mathcal{E})$ verifying
 \begin{enumerate}
\item $\theta_0=\hbox{Id}$, $\theta_t\circ \theta_s=
\theta_{t+s}$, for all $t,s\in \reali$
\item $\theta_t$ is measurable:
 $(\Omega,\mathcal{F}_t)\rightarrow
  (\Omega,\mathcal{F}_0)$ and $\{\{\theta_t\in A\}
  :A\in \mathcal{F}_0 \}=\mathcal{F}_t$
\item $\mathbb{P}\{\theta_t\in A\}= \mathbb{P}(A) $ for all $A\in
\mathcal{F}_0 $
\item $W_t\circ \theta_s=W_{t+s}-W_s$
\end{enumerate}
According to this framework we introduce the definition of
stationary stochastic process.
\begin{definition}\label{def-staz} We say that a stochastic process
$X: [0,\infty[\times \Omega \rightarrow \R^m$, is stationary if for all $s\in \reali$
$$X_t\circ \theta_s= X_{t+s}\qquad \hbox{$\mathbb{P}$-a.s. for a.e. $t\geq 0$}$$
\end{definition}
We assume all the
coefficients $A$, $B$, $C$, $D$ and $S$ to be
stationary stochastic processes. Namely on the coefficients we make the following assumptions:
\begin{hypothesis}
\label{genhyp} $\ $

\begin{enumerate}
\item [ A1)]$A:\left[  0,+\infty\right)  \times\Omega\rightarrow$
$\mathbb{R}^{n\times n}$, $B:\left[  0,+\infty\right)
\times\Omega \rightarrow$ $\mathbb{R}^{n\times k}$, $C^{i}:\left[
0,+\infty\right) \times\Omega\rightarrow$ $\mathbb{R}^{n\times
n}$, $i=1,...,d$ and $D^{i}:\left[  0,+\infty\right)
\times\Omega\rightarrow$ $\mathbb{R}^{n\times k}$, $i=1,...,d$,
are uniformly bounded process adapted to the filtration $\left\{
\mathcal{F}_{t}\right\}  _{t\geq0}$.

\item[ A2)] $S:\left[  0,+\infty\right)  \times\Omega\rightarrow
\mathbb{R}^{n\times n}$ is uniformly bounded and adapted to the
filtration $\left\{  \mathcal{F}_{t}\right\}  _{t\geq0}$ and it is
almost surely and almost everywhere symmetric and nonnegative. Moreover we assume that there exists $\beta >0$ such that $S\geq \beta I$.

\item[ A3)] $A$, $B$, $C$, $D$ and $S$ are stationary processes.
\end{enumerate}
\end{hypothesis}

In this case we immediately get:
\begin{lemma} Fix $T>0$ and let hypothesis \ref{genhyp} holds true.
Let $(P,Q)$ be the solution of the finite horizon BSRE
\begin{equation}
\left\{
\begin{array}
[c]{l} -dP_{t}=G\left(
A_{t},B_{t},C_{t},D_{t};S_{t};P_{t},Q_{t}\right)  dt+
{\displaystyle\sum_{i=1}^{d}}

Q_{t}^{i}dW_{t}^{i}, \qquad t \in [0,T]\\
P_{T}=P_T.
\end{array}
\right.\label{Riccati0T.staz}
\end{equation}
For fixed $s>0$ we define $\widehat{P}(t+s)=P(t)\theta_s$,
$\widehat{Q}(t+s)=Q(t)\theta_s$ then $(\widehat{P},\widehat{Q})$
is the unique solution in $[s,T+s]$ of the equation
\begin{equation}
\left\{
\begin{array}
[c]{l} -d\widehat{P}_{t}=G\left(
A_{t},B_{t},C_{t},D_{t};S_{t};\widehat{P}_{t},\widehat{Q}_{t}\right)
dt+ {\displaystyle\sum_{i=1}^{d}}

\widehat{Q}_{t}^{i}dW_{t}^{i}, \qquad t \in [s,T+s]\\
\widehat{P}_{T}=P_T\circ\theta_s.
\end{array}
\right.\label{Riccati0T.staz1}
\end{equation}
\end{lemma}
In the stationary assumptions the backward stochastic Riccati equation
\begin{align}
\label{RiccatiSTAZ}
&dP_{t}    =-\left[  A_{t}^{\ast}P_{t}+P_{t}A_{t}+S_{t}+
{\displaystyle\sum_{i=1}^{d}} \left(  \left(  C_{t}^{i}\right)
^{\ast}P_{t}C_{t}^{i}+\left(  C_{t} ^{i}\right)
^{\ast}Q_{t}+Q_{t}C_{t}^{i}\right)  \right]  dt+
{\displaystyle\sum_{i=1}^{d}}
Q_{t}^{i}dW_{t}^{i}+ \\ \nonumber
&  \left[  P_tB_t+ {\displaystyle\sum_{i=1}^{d}} \left(  \left(
C^{i}_t\right)  ^{\ast}P_t D^{i}_t +Q^{i}D_t ^{i}\right)  \right]
\left[  I+ {\displaystyle\sum_{i=1}^{d}} \left(  D^{i}_t\right)
^{\ast}P_tD^{i}_t\right]^{-1}\!\!\left[  P_tB_t+
{\displaystyle\sum_{i=1}^{d}} \left(  \left(  C^{i}_t\right)
^{\ast}P_t D_t ^{i}+Q^{i}_t D^{i}_t\right)
\right] ^{\ast}\!\!\,\!\! dt, \nonumber \\
\end{align}
admits a minimal solution $(\overline{P},\overline{Q})$, in the sense that
whenever another couple $(P,Q)$ is a solution to the Riccati equation then
$P -\overline{P}$ is a non-negative matrix, see also Corollary 3.3 in
\cite{GTinf} and definiton 3.2 in \cite{GM}.
This minimal solution $(\overline{P},\overline{Q})$ turns out to be stationary.
\begin{proposition}\label{prop-staz}
Assume hypothesis \ref{genhyp}, then the minimal solution
$(\overline{P},\overline{Q})$ of the infinite horizon stochastic
Riccati equation \eqref{RiccatiSTAZ} is stationary.
\end{proposition}
\Dim For all $\rho>0$ we denote by $P^{\rho}$ the
solution of equation (\ref{Riccati0T.staz})  in $[0,\rho]$ with final condition
$P^{\rho}(\rho)=0$. Denoting by $\lfloor \rho\rfloor$ the integer part of $\rho$,
we have, following Proposition 3.2 in \cite{GTinf} 
that for all $N$ for all $t\in [0,\lfloor N+s\rfloor]$, $P^{\lfloor
N+s\rfloor}_t\leq P^{ N+s}_t\leq P^{\lfloor N+s\rfloor+1}_t$,
$\mathbb{P}$-a.s.. Thus we can conclude noticing that by lemma
\ref{Riccati0T.staz1}
$$P^{N+s}_{t+s}=P^{N}_t\circ \theta_s.$$
Thus letting $N\rightarrow +\infty$ we obtain that for all $t\geq
0$, and $s>0$:
$$\mathbb{P}\left\{\overline{P}_{t+s}=\overline{P}_t\circ
\theta_s\right\}=1.$$ Now  $\overline{P}_{T+s}= \overline{P}_T
\circ \theta_s= \overline{P}_T$ so if one consider
\eqref{Riccati0T.staz} in the intervall $[s,T+s]$ with final data
$\overline{P}_{T+s}$ and \eqref{Riccati0T.staz1} with final data
$\overline{P}_T\circ \theta_s$, by the uniqueness of the solution
it follows that $Q_r=\hat{Q}_r, \ \mathbb{P}-\text{ a.s. and for almost
all } r \in [s, T+s]$.
 \finedim

We notice that in the BSRDE \eqref{RiccatiSTAZ} the final condition has been replaced by the stationarity condition on the solution process $(P,Q)$.

Next we give some definitions.
\begin{definition}
\label{defstab}We say that $(A,B,C,D)$ is stabilizable relatively to the
observations $\sqrt{S}$ (or $\sqrt{S}$-stabilizable) if there exists a control
$u\in L_{\mathcal{P}}^{2}([0,+\infty)\times\Omega;\R^k)$ such that for all
$t\geq0$ and all $x\in\mathbb{R}^{n}$%
\begin{equation}
\mathbb{E}^{\mathcal{F}_{t}}\int_{t}^{+\infty}[\left\langle S_{s}X_{s}%
^{t,x,u},X_{s}^{t,x,u}\right\rangle +|u_{s}|^{2}]ds<M_{t,x}.\label{condstabS}%
\end{equation}
for some positive constant $M_{t,x}$ where $X^{t,x,u}$ is the solution of the linear equation
\begin{equation}
\left\{
\begin{array}
[c]{ll}
dX_{s}=(A_{s}X_{s}+B_{s}u_{s})ds+
{\displaystyle\sum_{i=1}^{d}}
\left(  C_{s}^{i}X_{s}+D_{s}^{i}u_{s}\right)  dW_{s}^{i} & s\geq 0\\
X_{0}=x. &
\end{array}
\right.  \label{linear.stato.stato.intro}%
\end{equation}
\end{definition}
\noindent This kind of stabilizability condition, also called
finite cost condition, has been introduced in \cite{GTinf}. This
condition has been proved to be equivalent to the existence of a
minimal solution $(\bar P, \bar Q)$ of the Riccati equation
\eqref{RiccatiSTAZ}. Moreover whenever the first component $\bar
P$ is uniformly bounded in time it follows that the constant
$M_{t,x}$ appearing in \eqref{condstabS} can be chosen independent
of time.
\begin{definition}\label{def.stab.identita}
Let $P$ be a solution to equation \eqref{RiccatiSTAZ}. We say that
$P$ stabilizes $(A,B,C,D)$ relatively to the identity $I$ if for
every $t > 0$ and $x \in \R^n$ there exists a positive constant
$M$, independent of $t$, such that
\begin{equation}\label{condstabiden}
\Et\int_t^{+\infty}|X^{t,x}(r)|^2 dr \leq M \quad \quad
\mathbb{P}-a.s.,\end{equation} where $X^{t,x}$ is a mild
solution to:
\begin{equation}
\left\{
\begin{array}
[c]{ll}
d\overline{X}_{t}=\left[  A\overline{X}_{t}-B_{t}\left(  I+%

{\displaystyle\sum_{i=1}^{d}}

\left(  D_{t}^{i}\right)  ^{\ast}P_{t}D_{t}^{i}\right)  ^{-1}\left(
P_{t}B_{t}+

{\displaystyle\sum_{i=1}^{d}}
\left(  Q_{t}^{i}D_{t}^{i}+\left(  C_{t}^{i}\right)  ^{\ast}P_{t}D_{t}%
^{i}\right)  \right)  ^{\ast}\overline{X}_{t}\right]  dt +& \\
{\displaystyle\sum_{i=1}^{d}}\left[  C_{s}^{i}\overline{X}_{t}-D_{s}^{i}\left(  I+%

{\displaystyle\sum_{i=1}^{d}}
\left(  D_{t}^{i}\right)  ^{\ast}P_{t}D_{t}^{i}\right)  ^{-1}\left(
P_{t}B_{t}+
{\displaystyle\sum_{i=1}^{d}}

\left(  Q_{t}^{i}D_{t}^{i}+\left(  C_{t}^{i}\right)  ^{\ast}P_{t}D_{t}%
^{i}\right)  \right)  ^{\ast}\overline{X}_{t}\right]   dW_{t}, &
\\
\overline{X}_{0}=x &
\end{array}
\right.
\end{equation}
\end{definition}
From now on we assume that
\begin{hypothesis}\label{ipotesi2}
\item[(i)]  $(A,B,C,D)$ is $\sqrt{S}$-
stabilizable;
\item[(ii)] the process $\overline{P}$ is uniformly bounded in time;
\item[(iii)] the minimal solution $\bar{P}$ stabilizes $(A,B,C,D)$ with respect to the identity $I$.
\end{hypothesis}
We refer to \cite{GM} for cases when $P$
stabilizes $(A,B,C,D)$ relatively to the identity $I$. Notice that, thanks to the stationarity assumptions the
stabilizability condition can be simplified, see Remark 5.7 of
\cite{GTinf}.

Next we study the dual (costate) equation in the stationary case.
We denote by
\begin{align}
& \Lambda\left( t, \overline{P}_{t}, \overline{Q}_{t} \right)=-\left(  I+{\displaystyle\sum_{i=1}^{d}}\left(  D_{t}^{i}\right)  ^{\ast}\overline{P}_{t}D_{t}^{i}\right)  ^{-1}\left(\overline{P}_{t}B_{t}+{\displaystyle\sum_{i=1}^{d}}\left( \overline{Q}_{t}^{i}D_{t}^{i}+\left(  C_{t}^{i}\right)  ^{\ast}\overline{P}_{t}D_{t}^{i}\right)  \right)  ^{\ast}, \nonumber\\
& H_{t}=A_{t}+B_{t}\Lambda\left( t, \overline{P}_{t},\overline{Q}_{t} \right), \nonumber\\
& K_{t}^{i}=C_{t}^{i}+D_{t}^{i}\Lambda\left( t, \overline{P}_{t},\overline{Q}_{t} \right). \label{notazionifHK}
\end{align}
Thanks to Proposition \ref{prop-staz}, all the coefficients that appear in equation
\begin{equation}
\left\lbrace
\begin{array}[c]{ll}
dr_{t}=-H_{t}^{*}r_tdt-\bar{P}_{t}f_{t}dt-{\displaystyle\sum_{i=1}^{d}\left( K_{t}^{i} \right)^{*}g_{t}^{i}}dt+{\displaystyle\sum_{i=1}^{d}}g_{t}^{i}dW_{t}^{i}, & t \in \left[ 0,T\right]  \\
r_{T}=0. &
\end{array}
\right. \label{dualeT}
\end{equation}
are stationary so exactly as before we deduce that
for the solution $(r_T,g_T)$ the following holds:
\begin{lemma} Let $A$, $B$, $C$, $D$ and $S$ satisfy hypothesis \ref{genhyp} and let $f\in L^{\infty}_\P \left(\Omega\times\left[0,+\infty\right)\right)$ be a stationary process. Fix $T>0$ and $r_T \in L^\infty_\P(\Omega, \F_T;\R^n)$.
Let $(r,g)$ a solution to equation
\begin{equation}
\left\lbrace
\begin{array}[c]{ll}
dr_{t}=-H_{t}^{*}r_tdt-\bar{P}_{t}f_{t}dt-{\displaystyle\sum_{i=1}^{d}
\left( K_{t}^{i} \right)^{*}g_{t}^{i}}dt+{\displaystyle\sum_{i=1}^{d}}g_{t}^{i}dW_{t}^{i},
 & t \in \left[ 0,T\right]  \\
r_{T}=r_T. &
\end{array}
\right.
\end{equation}
For fixed $s>0$ we define $\widehat{r}_{t+s}=r_t\circ\theta_s$,
$\widehat{g}_{t+s}=g_t\circ\theta_s$ then $(\widehat{r},\widehat{g})$
is the unique solution in $[s,T+s]$ of the equation
\begin{equation}
\left\lbrace
\begin{array}[c]{ll}
d\widehat{r}_{t}=-H_{t}^{*}\widehat{r}_tdt-\bar{P}_{t}f_{t}dt-{\displaystyle\sum_{i=1}^{d}
\left( K_{t}^{i}
\right)^{*}\widehat{g}_{t}^{i}}dt+{\displaystyle\sum_{i=1}^{d}}\widehat{g}_{t}^{i}dW_{t}^{i},
 & t \in \left[ s,T+s\right]  \\
\widehat{r}_{T}=r _T \circ \theta_s. &
\end{array}
\right.
\end{equation}
\end{lemma}

Hence arguing as for the first component $\overline{P}$, we get that
the solution of the infinite horizon dual equation is stationary, as stated in the following proposition:
\begin{proposition}
\label{prop-staz2} Assume hypothesis \ref{genhyp} and hypothesis \ref{ipotesi2}, then the solution $(r^{\natural},g^{\natural})$ of
\begin{equation}
dr_{t}=-H_{t}^{*}r_t
dt-\overline{P}_{t}f_{t}dt-{\displaystyle\sum_{i=1}^{d}\left( K_{t}^{i}
\right)^{*}g_{t}^{i}}dt+{\displaystyle\sum_{i=1}^{d}}g_{t}^{i}dW_{t}^{i},
\label{dualeSTAZ}
\end{equation}
obtained as the
pointwise limit of the solution to equation \eqref{dualeT} is stationary. Moreover
$\left( r^\natural,g^\natural\right) \in L^{\infty}_\P\left( \Omega
\times \left[ 0,1\right],\R^{n}\right)\times L^2_\P\left(\Omega \times\left[ 0,1\right],
\R^{n\times d}\right)$.
\end{proposition}
\Dim The proof follows from an argument similar to the one in
Proposition 4.5 in \cite{GM}. Stationarity of the solution
$(r^{\natural}, g^{\natural})$ follows from the previous lemma.
 \finedim

Again we notice that in the dual BSDE \eqref{dualeSTAZ} the final condition has been replaced by the stationarity condition on the solution process $(r^\natural,g^\natural)$.

\noindent We need to show that in the stationary assumptions, the solution
of the closed loop equation is stationary. By using notation \eqref{notazionifHK},
we consider the following stochastic differential equation, which will
turn out to be the closed loop equation:
\begin{equation}
dX_{s}=H_{s}X_sds+{\displaystyle\sum_{i=1}^{d}}K_{s}^{i}X_{s}dW_{s}^{i}+B_s(B^*_sr^{\natural}_s
+{\displaystyle\sum_{i=1}^{d}}D_{s}^{i}g^{\natural,i}_{s})ds+f_sds+ {\displaystyle\sum_{i=1}^{d}}D^i_s (B^*_sr^{\natural}_s+{\displaystyle\sum_{i=1}^{d}}D_{s}^{i}g^{\natural,i}_{s})dW^i_s,
\label{closedloopSTAZ}
\end{equation}
where $(r^{\natural},g^{\natural})$ is the solution of the dual (costate) equation \eqref{dualeSTAZ}.

\begin{proposition}
\label{propclosedloopSTAZ}
 Assume hypothesis \ref{genhyp}  and hypothesis \ref{ipotesi2} holds true
then there exists a unique stationary solution of equation \eqref{closedloopSTAZ}.
\end{proposition}
\Dim We set $f^1_s=f_s+B_s(B^*_sr^{\natural}_s
+{\displaystyle\sum_{i=1}^{d}}D_{s}^{i}g^{\natural,i}_{s})$ and
$f^{2,j}_s=D^j_s(B^*_sr^{\natural}_s
+{\displaystyle\sum_{i=1}^{d}}D_{s}^{i}g^{\natural,i}_{s})$,
$j=1,...,d$. We can extend $f^1$, $f^2$ for negative times letting
for all $t\in[0,1]$, $f^i_{-N+t}=f^i_{t}\circ \theta_{-N}$,
$i=1,2$, $N\in \nat$. We notice that $f^i\mid_{[-N, +\infty)}$ is
predictable with respect to the filtration $(\G^{-N}_t)_{t\geq
-N}$. Therefore for all $N\in \nat$ equation
\begin{equation*}
\left\{
\begin{array}
[c]{l}
dX^{-N}_{s}=H_{s}X^{-N}_{s}ds+\ +\sum_{i=1}^{d}  K_{s}^{i}X^{-N}_{s}dW_{s}^{i}
+f^1_s ds+{\displaystyle\sum_{i=1}^{d}}f^{2,i}_sdW_{s}^{i},\\
X^{-N}_{-N}=0,
\end{array}
\right.
\end{equation*}
admits a solution $(X^{-N,0}_t)_t$, defined for $t\geq -N$ and predictable
with respect to the filtration $(\G^{-N}_t)_{t\geq -N}$. We extend $X^{-N,0}$
to the whole real axis by setting $X^{-N,0}_t=0$ for $t<-N$. We want to
prove that, fixed $t\in \reali$, $(X^{-N,0}_t)_N$ is a Cauchy sequence
in $L^2(\Omega)$. In order to do this we notice that for $t\geq -N+1$, $X^{-N,0}_t-X^{-N+1,0}_t$
solves the following (linear) stochastic differential equation
\begin{equation*}
X^{-N,0}_t-X^{-N+1,0}_t=X^{-N,0}_{-N+1}+\int_{-N+1}^t H_s (X^{-N,0}_s-X^{-N+1,0}_s)
+\displaystyle\sum_{i=1}^d\int_{-N+1}^t K^i_s(X^{-N,0}_s-X^{-N+1,0}_s)dW^i_s.
\end{equation*}
By the Datko theorem, see e.g. \cite{GM} and \cite{GTinf}, there exist constants $a,c>0$ such that
\begin{equation*}
(\E\vert X^{-N,0}_t-X^{-N+1,0}_t\vert^2)^{1/2}\leq C e^{-\frac{a(t+N-1)}{2}}(\E\vert X^{-N,0}_{-N+1}\vert^2)^{1/2}.
\end{equation*}
So, fixed $t\in \reali$ and $M,N \in \nat$, $M>N$ sufficiently large such that $-N\leq t$,
\begin{equation}
(\E\vert X^{-N,0}_t-X^{-M,0}_t\vert^2)^{1/2}\leq  \sum _{k=N}^{M-1}(\E\vert X^{-k,0}_t-X^{-k+1,0}_t\vert^2)^{1/2}
\leq  C \sum _{k=N}^{M-1}e^{-\frac{a(t+k-1)}{2}}(\E\vert X^{-k,0}_{-k+1}\vert^2)^{1/2}.\label{stima1}
\end{equation}
Next we look for a uniform estimate with respect to $k$ of $\E\vert X^{-k,0}_{-k+1}\vert^2$. For $s\in[-k,-k+1]$,
\begin{equation}
\label{eq1}
X^{-k,0}_s=\int_{-k}^sA_rX^{-k,0}_rdr+\int_{-k}^sB_r\bar{u}_rdr+
\sum_{i=1}^d\int_{-k}^sC^i_rX^{-k,0}_rdW^i_r+\sum_{i=1}^d\int_{-k}^sD^i_r\bar{u}_rdW^i_r+\int_{-k}^sf_rdr,
\end{equation}
where $\bar{u}$ is the optimal control that minimizes the cost
\begin{equation*}
J(-k,0,u)=\E\int_{-k}^{-k+1}[\vert \sqrt{S_s}X_s \vert^2+\vert u_s\vert^2] ds.
\end{equation*}
By computing $d[\langle \overline{P}_s X^{-k,0}_s, X^{-k,0}_s\rangle
+2\langle \bar r^{\natural} _s,X^{-k,0}_s\rangle]$ we get, for every $T>0$,
\begin{align*}
&\E\int_{-k}^{-k+1} [\vert \sqrt{S_s}X_s \vert^2 +|\bar{u}_{s}|^{2}]ds
=-\E \langle \overline{P}_{-k+1} X^{-k,0}_{-k+1},X^{-k,0}_{-k+1}\rangle-2\E \int_{-k}^{-k+1} \langle r^{\natural} _s,f_s\rangle ds \nonumber \\
&-\E \int_{-k}^{-k+1} \vert\left(  I+\sum_{i=1}^{d}\left(  D_{s}^{i}\right)^*\overline{P}_s D_s^i \right) ^{-1}( B^*_s r^{\natural}_s+\sum_{i=1}^{d}\left(  D_{s}^{i}\right)  ^{\ast} g_{s}^{\natural,i})\vert  ^2 ds
 \leq 2\vert \E \int_{-k}^{-k+1} \langle r^{\natural}_s,f_s\rangle ds\vert \leq A, \nonumber \\
\end{align*}
where $A$ is a constant independent on $k$.
By \eqref{eq1} we get
\begin{equation*}
\sup _{-k\leq s\leq -k+1}\E \vert X^{-k,0}_s \vert^2 \leq
C\int_{-k}^s \sup _{-k\leq r\leq s}\E \vert X^{-k,0}_r \vert^2
dr+C\E \int_{-k}^{-k+1} \vert \bar{u}_r \vert ^2dr+\E
\int_{-k}^{-k+1} \vert f_r \vert^2 dr,
\end{equation*}
and so by applying the Gronwall lemma, we get
\begin{equation*}
\sup _{-k\leq s\leq -k+1}\E \vert X^{-k,0}_s \vert^2
\leq C e^C(A+\E \int_{-k}^{-k+1} \vert f_r \vert^2 dr).
\end{equation*}
Since $f$ is stationary, we can conclude that
\begin{equation*}
\sup _{-k\leq s\leq -k+1}\E \vert X^{-k,0}_s \vert^2 \leq C,
\end{equation*}
where $C$ is a constant independent on $k$. By \eqref{stima1},
we get
\begin{equation*}
(\E\vert X^{-N,0}_t-X^{-M,0}_t\vert^2)^{1/2}
\leq C\sum _{k=N}^{M-1}e^{-\frac{a(t+k-1)}{2}}.
\end{equation*}
So we can conclude that, fixed $t \in \reali$, $(X^{-N,0}_t)_N$ is a
Cauchy sequence in $L^2(\Omega)$, and so it converges in $L^2(\Omega)$ to a random
variable denoted by $\zeta^{\natural}_t$. Notice that for every $t \in \reali$
we can define $\zeta^{\natural}_t$, and we prove that $\zeta^{\natural}$ is a
stationary process.  Let $t\in \reali$, $-N<t$ and $s>0$:
since the shift $\theta $ is measure preserving,
\begin{equation*}
\lim_{N\rightarrow \infty}\E\vert X^{-N,0}_t\circ \theta_s-\zeta^{\natural}_t\circ \theta_s\vert^2=0,
\end{equation*}
moreover $X^{-N,0}_t\circ \theta_s=X^{-N+s,0}_{t+s}$ and
\begin{equation*}
\lim_{N\rightarrow \infty}\E\vert X^{-N+s,0}_{t+s}-\zeta^{\natural}_{t+s}\vert^2=0.
\end{equation*}
By uniqueness of the limit we conclude that $\zeta^{\natural}_t\circ \theta_s=\zeta^{\natural}_{t+s}$.
Notice that since $N\in \nat$ and $\F_0\supset \G^{-N}_{0}$, $\zeta^{\natural}_0$ is
$\F_0$-measurable. Let us consider the value of the solution of equation
\eqref{closedloopSTAZ} starting from $X_0=\zeta^{\natural}_0$. By stationarity
of the coefficients and of $\zeta^{\natural}$, we get that $X$ is a stationary
solution of equation \eqref{closedloopSTAZ}, that we denote by $X^{\natural}$.
In order to show the uniqueness of the periodic solution it is enough to
notice that if $f^j=0$, $j=1,2$, and $X^{\natural}$ is a periodic
solution of \eqref{closedloopSTAZ}, then
\begin{equation*}
\E\vert X^{\natural}_{0}\vert^2=\E\vert X^{\natural}_{N}\vert^2
\leq Ce^{-aN}\E\vert \zeta^{\natural}_0\vert^2.
\end{equation*}
Therefore $X^{\natural}_0=0$ and this concludes the proof.
\finedim

We can now treat the following optimal control problem for a stationary cost functional: minimize
over all admissible controls $u\in\U^{\natural}$ the cost functional
\begin{equation}
J^{\natural}(u,X)=\E\int_{0}^{1}[\vert \sqrt{S_s}X_s \vert^2+\vert u_s\vert^2] ds, \quad (u,X)\in \U^{\natural},
\label{costoSTAZ}
\end{equation}
where
\begin{equation}
\U^{\natural}=\left\lbrace (u,X)\in L^2_\P(\Omega\times[0,1])\times C([0,1],L^2_\P(\Omega)):
X_s=X_0 \circ \theta_s, \text{ \\\\} \forall s\in \reali
\right\rbrace
\end{equation}
and $X$ is the solution of equation
\begin{equation}
dX_t=A_tX_tdt+B_tu_tdt+\sum_{i=1}^dC^i_tX_tdW^i_t+\sum_{i=1}^dD^i_tu_tdW^i_t+f_tdt,
\end{equation}
relative to $u$.

\begin{theorem}
\label{teocontrolloSTAZ}
Let $X^{\natural}\in C([0,1],L^2(_\P\Omega))$ be the unique stationary solution of equation \eqref{closedloopSTAZ} and let
\begin{equation}
u^{\natural}_{t}=\!-\!\left(  I+ \sum_{i=1}^{d} \left(
D_{t}^{i}\right)  ^{\ast}\overline{P}_{t}D_{t}^{i}\right)
^{-1}\!\!\!\left(\overline{ P}_{t}B_{t}+
\sum_{i=1}^{d}\left(\overline{Q}_{t}^{i}D_{t}^{i}+\left(  C_{t}^{i}\right)  ^{\ast}\overline{Q}_{t}D_{t}%
^{i}\right)  \right)  ^{\ast}\!\!\!X^{\natural}_{t}+B_{t}^{*}r^{\natural}_{t}+\! \! \sum_{i=1}^{d} (D_{t}^{i})^*g^{\natural ,i}_{t}. \label{feedback.STAZ}%
\end{equation}
Then $(u^{\natural},X^{\natural})\in\U^{\natural}$ and it is the unique optimal couple for the cost \eqref{costoSTAZ}, that is
\begin{equation*}
J^{\natural}(u^{\natural},X^{\natural})=\inf_{(u,X)\in \U^{\natural}}J^{\natural}(u,X).
\end{equation*}
The optimal cost is given by
\begin{equation}
\label{costoottimoSTAZ}
\overline{J}^{\natural}=J^{\natural}(u^{\natural},X^{\natural})
=2\E\int_0 ^{1} \langle r^{\natural}_s ,f_s\rangle ds
 -\E\int_0 ^{1} \vert (  I+{\displaystyle\sum_{i=1}^{d}}
\left(  D_{t}^{i}\right)  ^{\ast}\overline{P}_{t}D_{t}^{i})  ^{-1}( B_{t}^{*}r^{\natural}_{t}+{\displaystyle\sum_{i=1}^{d}}\left( D_{t}^{i}\right) ^{*}g_{t}^{\natural,i})\vert ^2 ds.
\end{equation}

\end{theorem}

\Dim By computing $d\langle \overline{P}_s X_s, X_s\rangle +2\langle r^{\natural} _s,X_s\rangle$ we get
\begin{align*}
& \E\int_0^{1} [\left\langle S_{s}X_{s},X_{s}\right\rangle
+|u_{s}|^{2}]ds = \E\langle\overline{P}_0 X_0,X_0\rangle-\E
\langle \overline{P}_1 X_1, X_1\rangle +2\E\langle r^{\natural}_0
,X_0\rangle
-2\E\langle r^{\natural}_1 ,X_1\rangle-2\E \int_0^1 \langle  r^{\natural}_s,f_s\rangle ds \nonumber \\
& +\E \int_0^1 \vert \left(  I+\sum_{i=1}^{d} \left(
D_{s}^{i}\right)  ^{\ast}\overline{P}_{s} D_{s}^{i}\right)
^{1/2}\left(u_s +( I+\sum_{i=1}^{d}
\left(  D_{s}^{i}\right)  ^{\ast}\overline{P}_{s} D_{s}^{i})  ^{-1}\right.* \nonumber \\
&\left. *\left(
\overline{P}_{s} B_{s}+{\displaystyle\sum_{i=1}^{d}}\left(\overline{Q}_{s}^{i}D_{s}^{i}+\left(  C_{s}^{i}\right)  ^{\ast}\overline{P}_{s}D_{s}^{i}\right)  \right)  ^{\ast}{X}_{s}
  +B_s ^* r^{\natural} _s +\sum_{i=1}^d D^i_s(\bar g^{\natural,i}_s)^*\right)\vert ^2 ds \nonumber \\
&-\E \int_0^1 \vert\left(  I+\sum_{i=1}^{d}\left(  D_{s}^{i}\right)^*\overline{P}_s D_s^i \right) ^{-1}( B^*_s r^{\natural}_s+\sum_{i=1}^{d}\left(  D_{s}^{i}\right)  ^{\ast}\bar g_{s}^{\natural,i})\vert  ^2 ds.\nonumber \\
\end{align*}
Since by Propositions \ref{prop-staz}, \ref{prop-staz2}, and \ref{propclosedloopSTAZ} $(u,X)\in \U^{\natural}$,
we get
\begin{align*}
& \E\int_0^{1} [\left\langle S_{s}X_{s},X_{s}\right\rangle +|u_{s}|^{2}]ds =-2\E \int_0^1 \langle  r^{\natural}_s,f_s\rangle ds  +\E \int_0^1 \vert \left(  I+\sum_{i=1}^{d}
\left(  D_{s}^{i}\right)  ^{\ast}\overline{P}_{s} D_{s}^{i}\right)  ^{1/2}\times \nonumber \\
&\times\left(u_s +( I+\sum_{i=1}^{d}
\left(  D_{s}^{i}\right)  ^{\ast}\overline{P}_{s} D_{s}^{i})  ^{-1}
\left(
\overline{P}_{s} B_{s}+{\displaystyle\sum_{i=1}^{d}}\left(\overline{Q}_{s}^{i}D_{s}^{i}+\left(  C_{s}^{i}\right)  ^{\ast}\overline{P}_{s}D_{s}^{i}\right)  \right)  ^{\ast}{X}_{s}
  +B_s ^* r^{\natural} _s +\sum_{i=1}^d D^i_s(\bar g^{\natural,i}_s)^*\right)\vert ^2 ds \nonumber \\
&-\E \int_0^1 \vert\left(  I+\sum_{i=1}^{d}\left(  D_{s}^{i}\right)^*\overline{P}_s D_s^i \right) ^{-1}( B^*_s r^{\natural}_s+\sum_{i=1}^{d}\left(  D_{s}^{i}\right)  ^{\ast}\bar g_{s}^{\natural,i})\vert  ^2 ds.\nonumber \\
\end{align*}
So
\begin{equation}
u^{\natural}_{t}=\!-\!\left(  I+ \sum_{i=1}^{d} \left(
D_{t}^{i}\right)  ^{\ast}\overline{P}_{t}D_{t}^{i}\right)
^{-1}\!\!\!\left( \overline{P}_{t}B_{t}+
\sum_{i=1}^{d}\left(\overline{Q}_{t}^{i}D_{t}^{i}+\left(  C_{t}^{i}\right)  ^{\ast}\overline{Q}_{t}D_{t}%
^{i}\right)  \right)  ^{\ast}\!\!\!X^{\natural}_{t}+B_{t}^{*}r^{\natural}_{t}+\! \! \sum_{i=1}^{d} (D_{t}^{i})^*g^{\natural ,i}_{t}
\end{equation}
is the optimal cost: $u^{\natural}$ minimizes the cost \eqref{costoottimoSTAZ}, and the corresponding state $X^{\natural}$ is stationary by Proposition \ref{propclosedloopSTAZ}, so that $(u^{\natural},X^{\natural})\in\U^{\natural}$.
\finedim

\section{Ergodic control}
In this section we consider cost functionals depending only on the
asymptotic behaviour of the state (ergodic control).
Throughout this section we assume the following:
\begin{hypothesis}\label{ipotesi_erg}
The coefficient satisfy hypothesis \ref{genhyp}, and moreover
\begin{itemize}
\item $S\geq \epsilon I$, for some $\epsilon >0$. \item
$(A,B,C,D)$ is stabilizable relatively to $S$. \item The first
component of the minimal solution $P$ is bounded in
time.
\end{itemize}
\end{hypothesis}
Notice that these conditions implies that
$(P,Q)$ stabilize $(A,B,C,D)$ relatively to the identity.

We first consider discounted cost functional and then we
compute a suitable limit of the discounted cost.
Namely, we consider the discounted cost
functional
\begin{equation}
J_{\alpha}(0,x,u)=\mathbb{E}\int_{0}^{+\infty}e^{-2\alpha s}[\left\langle S_{s}X_{s}%
^{0,x,u},X_{s}^{0,x,u}\right\rangle +|u_{s}|^{2}]ds,
\label{costo_scontato}
\end{equation}
where $X$ is solution to equation
\begin{equation*}
\left\{
\begin{array}
[c]{ll}
dX_{s}=(A_{s}X_{s}+B_{s}u_{s})ds+
{\displaystyle\sum_{i=1}^{d}}
\left(  C_{s}^{i}X_{s}+D_{s}^{i}u_{s}\right)  dW_{s}^{i} +f_s ds & s\geq t\\
X_{t}=x. &
\end{array}
\right.
\end{equation*}
$A$, $B$, $C$ and $D$ satisfy hypothesis \ref{genhyp} and $f\in
L_\P^{\infty}(\Omega\times[0,+\infty))$ and is a stationary process. When the coefficients are
deterministic the problem has been extensively studied, see e.g.
\cite{Ben1} and \cite{Tess2}.

Our purpose is to minimize the discounted cost functional with respect to every
admissible control $u$. We define the set of admissible controls as
\begin{equation*}
\mathcal U ^{\alpha}=\left\lbrace u\in L^2(\Omega\times[0,+\infty))
:\mathbb{E}\int_{0}^{+\infty}e^{-2\alpha s}[\left\langle S_{s}X_{s}%
^{0,x,u},X_{s}^{0,x,u}\right\rangle +|u_{s}|^{2}]ds <
+\infty\right\rbrace .
\end{equation*}
Fixed $\alpha >0$, we define $X^{\alpha}_s=e^{-\alpha s}X_s$
and $u^{\alpha}_s=e^{-\alpha s}u_s$. Moreover
we set $A^{\alpha}_s=A_s-\alpha I$ and
$f^{\alpha}_s=e^{-\alpha s}f_s$, and $f^{\alpha}\in
L^2_\P(\Omega\times[0,+\infty))\cap L^{\infty}_\P(\Omega\times[0,+\infty))$. $X^{\alpha}_s$
is solution to equation
\begin{equation}
\left\{
\begin{array}
[c]{ll}
dX^{\alpha}_s=(A^{\alpha}_sX^{\alpha}_s+B_su^{\alpha}_s)ds+
{\displaystyle\sum_{i=1}^{d}}
\left(  C_s^{i}X^{\alpha}_s+D_s^{i}u^{\alpha}_s\right)  dW_{s}^{i} +f^{\alpha}_s ds & s\geq 0\\
X^{\alpha}_{0}=x, &
\end{array}
\right.  \label{stato.alfa}%
\end{equation}
By the definition of $X^{\alpha}$, we note that if $(A,B,C,D)$ is stabilizable with respect to the identity,
then $(A^{\alpha},B,C,D)$ also is.
We also denote by $(P^{\alpha},Q^{\alpha})$ the minimal solution
of a stationary backward Riccati equation \eqref{RiccatiSTAZ}
with $A^{\alpha}$ in the place of $A$.
Since, for $0<\alpha<1$, $A^{\alpha}$ is uniformly bounded
in $\alpha$, also $P^{\alpha}$ is uniformly bounded
in $\alpha$. Arguing as in Proposition \ref{prop-staz},
$(P^{\alpha},Q^{\alpha})$ is a stationary process.

Let us denote by $(r^{\alpha},g^{\alpha})$ the solution of the infinite horizon BSDE
\begin{equation}
dr^{\alpha}_{t}=-(H^{\alpha}_{t})^{*}r^{\alpha}_t dt-P^{\alpha}_{t}f^{\alpha}_{t}dt-{\displaystyle\sum_{i=1}^{d}
\left( K^{\alpha,i}_{t} \right)^{*}g_{t}^{\alpha,i}}dt+{\displaystyle\sum_{i=1}^{d}}g_{t}^{\alpha,i}dW_{t}^{i},\text{
\ \ \ \ \ } t \geq 0 , \\
\label{dualeERG}
\end{equation}
where $H^{\alpha}$ and $K^{\alpha}$ are defined as in (\ref{notazionifHK}),
with $A^{\alpha}$, $P^{\alpha}$ and $Q^{\alpha}$ respectively in the place of $A$, $P$ and $Q$.
By \cite{GM}, section 4, we get that equation \eqref{dualeERG} admits a
solution $(r^{\alpha},g^{\alpha})\in L^2_\P(\Omega \times[0,+\infty))
\cap L^2_\P(\Omega \times[0,+\infty))
\times L^{\infty}_\P(\Omega \times[0,T])$, for every fixed $T>0$.

Moreover by \cite{GM}, section 6, we know that
\begin{align*}
&\underline{\lim}_{\alpha\rightarrow 0}\alpha\!\inf
_{u^{\alpha}\in\U^{\alpha}}J_{\alpha}\left( 0,x,u^{\alpha}\right)=
\\ \nonumber &\underline{\lim}_{\alpha\rightarrow
0}[\alpha\int_0^{+\infty}2\E\langle r^{\alpha}_s, f^{\alpha}_s\rangle
ds -\alpha\E\int_0 ^{+\infty} \vert (I+{\displaystyle\sum_{i=1}^{d}}
\left(  D_{s}^{i}\right)  ^{\ast}P^{\alpha}_{s}D_{s}^{i})
^{-1}(B_{s}^{*}r^{\alpha}_{s}+{\displaystyle\sum_{i=1}^{d}}\left(
D_{s}^{i}\right) ^{*}g_{s}^{\alpha,i})\vert ^2 ds ].\\ \nonumber
\end{align*}

We can also prove the following convergence result for $(r^{\alpha},g^{\alpha})$.
\begin{lemma}
 For all fixed $T>0$, $r^{\alpha}\mid_{[0,T]}\rightarrow r^{\natural}\mid_{[0,T]}$
in $L^2_\P (\Omega \times [0,T])$. Moreover, for every fixed $T>0$, as $\alpha\rightarrow 0$:
\begin{multline*}
 \E\int_0 ^{T} \vert (  I+
{\displaystyle\sum_{i=1}^{d}} \left(  D_{s}^{i}\right)
^{\ast}P^{\alpha}_{s}D_{s}^{i})  ^{-1}(
B_{s}^{*}r^{\alpha}_{s}+{\displaystyle\sum_{i=1}^{d}} \left(
D_{s}^{i}\right) ^{*}g_{s}^{\alpha,i})\vert ^2 ds \rightarrow\\
\E\int_0 ^{T} \vert (  I+ {\displaystyle\sum_{i=1}^{d}} \left(
D_{s}^{i}\right)  ^{\ast}\overline{P}_{s}D_{s}^{i})  ^{-1}(
B_{s}^{*}r^{\natural}_{s}+{\displaystyle\sum_{i=1}^{d}}\left(
D_{s}^{i}\right) ^{*}g_{s}^{\natural ,i})\vert ^2 ds
\end{multline*}

\end{lemma}
\Dim
The first assertion follows from lemma 6.6 in \cite{GM}.
Notice that stationarity of the coefficients in the limit equation gives
stationarity of the solution, and so it allows to identify the limit with
the stationary solution of the dual BSDE.
For the second assertion for the optimal couple $(X^\alpha,u^\alpha)$ for the
optimal control problem on the time interval $[0,T]$:
\begin{align}
&\int_0^T[\vert \sqrt{S_s}X^\alpha_s\vert^2+\vert u^\alpha_s\vert^2 ds
=\langle P^{\alpha}_{0}x,x\rangle+2\langle r^{\alpha}_0 ,x\rangle+
2\E\int_0 ^T \langle r^{\alpha}_s ,f^{\alpha}_s\rangle ds \nonumber \\
&\E\langle P^{\alpha}_{T}X^\alpha_T,X^\alpha_T\rangle+2\E\langle r^{\alpha}_T ,X^\alpha_T\rangle
-\E\int_0 ^T \vert (  I+
{\displaystyle\sum_{i=1}^{d}}
\left(  D_{t}^{i}\right)  ^{\ast}P^{\alpha}_{t}D_{t}^{i})  ^{-1}
( B_{t}^{*}r^{\alpha}_{t}+{\displaystyle\sum_{i=1}^{d}}\left( D_{t}^{i}\right)
^{*}g_{t}^{\alpha,i})\vert ^2 ds. \label{rel.fond.scontata}
\end{align}
Since, as $\alpha\rightarrow 0$, in (\ref{rel.fond.scontata}) all the
terms but the last one converge to the corresponding stationary term,
and since by \cite{GM} $(r^\alpha,g^\alpha)$ is uniformly,
with respect to $\alpha$, bounded
in $L^2_\P (\Omega \times [0,T])\times L^2_\P (\Omega \times [0,T])$, then
$(r^\alpha \mid_{[0,T]},g^\alpha \mid_{[0,T]})\rightharpoonup
(r^\natural\mid_{[0,T]},g^\natural\mid_{[0,T]})$ in
$L^2_\P (\Omega \times [0,T])\times L^2_\P (\Omega \times [0,T])$,
we get the desired convergence.

\finedim

\noindent This is enough to characterize the ergodic limit. Indeed we have that:
\begin{theorem} 
We get the following characterization of the optimal
cost:
\[
\lim_{\alpha\to 0} 2\alpha \inf_{u \in \mathcal{U}^{\alpha}} J_\alpha(x,u)= \E\left [2\langle f(0),r^{\natural}(0)\rangle
-\vert (I+{\displaystyle\sum_{i=1}^{d}}
\left(  D_{0}^{i}\right)  ^{\ast}\overline{P}_{0}D_{0}^{i})  ^{-1}(B_{0}^{*}r^{\natural}_{0}+{\displaystyle\sum_{i=1}^{d}}\left( D_{0}^{i}\right) ^{*}g_{0}^{\natural,i})\vert ^2  \right ].
\]
\end{theorem}
\Dim Let us define $\widetilde{r}^{\alpha}_{t}=e^{\alpha t}r^{\alpha}_{t}$, $\widetilde{g}^{\alpha}_{t}=e^{\alpha t}g^{\alpha}_{t}$. $(\widetilde{r}^{\alpha}_{t},\widetilde{g}^{\alpha}_{t})$ is the solution to
\begin{equation*}
d\widetilde{r}^{\alpha}_{t}=-(H^{\alpha}_{t})^{*}\widetilde{r}^{\alpha}_t dt+\alpha I\widetilde{r}^{\alpha}_{t}dt-P^{\alpha}_{t}f_{t}dt-{\displaystyle\sum_{i=1}^{d}\left( K^{\alpha,i}_{t} \right)^{*}\widetilde{g}_{t}^{\alpha,i}}dt+{\displaystyle\sum_{i=1}^{d}}\widetilde{g}_{t}^{\alpha,i}dW_{t}^{i},\text{
\ \ \ \ \ } t \geq 0 , \\
\end{equation*}
and so, arguing as in lemma \ref{prop-staz2}, $(\widetilde{r}^{\alpha}_{t},\widetilde{g}^{\alpha}_{t})$ are stationary processes.
Now we compute
\begin{align*}
\underline{\lim}_{\alpha\rightarrow 0}2\alpha\!\inf
_{u^{\alpha}\in\U^{\alpha}}J_{\alpha}\left( 0,x,u^{\alpha}\right)=
&\underline{\lim}_{\alpha\rightarrow
0}\left[ 2\alpha\int_0^{+\infty}e^{-2\alpha s}2\E\langle
\widetilde{r}^{\alpha}_s, f_s\rangle ds  \right. \\ \nonumber
&\left. -2\alpha\int_0 ^{+\infty} e^{-2\alpha s}\E\vert
(I+{\displaystyle\sum_{i=1}^{d}} \left(  D_{s}^{i}\right)
^{\ast}P^{\alpha}_{s}D_{s}^{i})
^{-1}(B_{s}^{*}\widetilde{r}^{\alpha}_{s}+{\displaystyle\sum_{i=1}^{d}}\left(
D_{s}^{i}\right) ^{*}\widetilde{g}_{s}^{\alpha,i})\vert ^2 ds \right]\\ \nonumber
& =\underline{\lim}_{\alpha\rightarrow
0}\left[ 2\alpha\sum_{k=1}^{\infty}e^{-2\alpha k}\int_0^1e^{-2\alpha
s}2\E\langle \widetilde{r}^{\alpha}_s, f_s\rangle ds \right. \\ \nonumber
&\left. -2\alpha\sum_{k=1}^{\infty}e^{-2\alpha k}\int_0 ^1 e^{-2\alpha
s}\E\vert (I+{\displaystyle\sum_{i=1}^{d}} \left( D_{s}^{i}\right)
^{\ast}P^{\alpha}_{s}D_{s}^{i})
^{-1}(B_{s}^{*}\widetilde{r}^{\alpha}_{s}+{\displaystyle\sum_{i=1}^{d}}\left(
D_{s}^{i}\right) ^{*}\widetilde{g}_{s}^{\alpha,i})\vert ^2 ds \right]\\ \nonumber
&= \underline{\lim}_{\alpha\rightarrow
0}\left[ 2\alpha\sum_{k=1}^{\infty}e^{-2\alpha k}\int_0^1 2\E\langle
r^{\alpha}_s, f^{\alpha}_s\rangle ds  \right. \\ \nonumber
&\left. -2\alpha\sum_{k=1}^{\infty}e^{-2\alpha k}\int_0 ^1\E \vert
(I+{\displaystyle\sum_{i=1}^{d}} \left(  D_{s}^{i}\right)
^{\ast}P^{\alpha}_{s}D_{s}^{i})
^{-1}(B_{s}^{*}r^{\alpha}_{s}+{\displaystyle\sum_{i=1}^{d}}\left(
D_{s}^{i}\right) ^{*}g_{s}^{\alpha,i})\vert ^2 ds \right] .\\ \nonumber
\end{align*}
Since $(r^{\alpha}_s, g^{\alpha}_s)\rightarrow (r^{\natural}_s, g^{\natural}_s)$
in $L^2_\P(\Omega \times[0,1])\times L^2_\P(\Omega \times[0,1])$ we get that
\begin{align*}
\underline{\lim}_{\alpha\rightarrow 0}2\alpha\!
\inf _{u^{\alpha}\in\U^{\alpha}}J_{\alpha}\left(  0,x,u^{\alpha}\right)
& =2\E\int_0^1\langle r^{\natural}_s, f_s\rangle ds
-\E\int_0 ^1 \vert (I+{\displaystyle\sum_{i=1}^{d}}
\left(  D_{s}^{i}\right)  ^{\ast}{P}^{\alpha}_{s}D_{s}^{i})  ^{-1}(B_{s}^{*}r^{\natural}_{s}+{\displaystyle\sum_{i=1}^{d}}
\left( D_{s}^{i}\right) ^{*}g_{s}^{\natural,i})\vert ^2 ds \\ \nonumber
&= 2\E\langle r^{\natural}_0, f_0\rangle-\E\vert (I+{\displaystyle\sum_{i=1}^{d}}
\left(  D_{0}^{i}\right)  ^{\ast}\overline P_{0}D_{0}^{i})  ^{-1}(B_{0}^{*}r^{\natural}_{0}+{\displaystyle\sum_{i=1}^{d}}
\left( D_{0}^{i}\right) ^{*}g_{0}^{\alpha,i})\vert ^2 ,\\ \nonumber
\end{align*}
where the first equality holds also in the periodic case and the second equality holds only in the stationary case.
\finedim

The next step is to minimize
\[
\widehat{J}(x,u)=\underline{\lim}_{\alpha\rightarrow 0} 2\alpha J(x,u)
\]
over all $u\in \widehat{\U}$, where
\begin{equation}
 \widehat{\U}=\left\lbrace u\in L^2_{loc}:\text{  \\}
\mathbb{E}\int_{0}^{+\infty}e^{-2\alpha s}[\left\langle S_{s}X_{s}%
^{0,x,u},X_{s}^{0,x,u}\right\rangle +|u_{s}|^{2}]ds<+\infty, \text{  }\forall \alpha >0.
 \right\rbrace
\label{Ucappuccio}
\end{equation}
We will prove that
\[
\inf_{u\in \widehat{\U}}\widehat{J}(x,u)=J^{\natural}(u).
\]
Let $\widehat{X}$ be solution of
\begin{equation*}
\left\lbrace
\begin{array}
[c]{l}
d\widehat{X}^x_{s}=H_{s}\widehat{X}^x_sds+{\displaystyle\sum_{i=1}^{d}}K_{s}^{i}\widehat{X}^x_{s}dW_{s}^{i}+B_s(B^*_sr^{\natural}_s
+{\displaystyle\sum_{i=1}^{d}}D_{s}^{i}g^{\natural,i}_{s})ds+f_sds+ {\displaystyle\sum_{i=1}^{d}}D^i_s (B^*_sr^{\natural}_s+{\displaystyle\sum_{i=1}^{d}}D_{s}^{i}g^{\natural,i}_{s})dW^i_s\\
\widehat{X}^x_0=x,
\end{array}
\right.
\end{equation*}
and let
\begin{equation*}
\widehat{u}^x_{s}=-\Lambda(s,\overline{P}_{s},\overline{Q}_s)\widehat{X}^x_{s}+(B^*_sr^{\natural}_s
+{\displaystyle\sum_{i=1}^{d}}D_{s}^{i}g^{\natural,i}_{s}).
\end{equation*}
Notice that by proposition \ref{propclosedloopSTAZ} if $x=\zeta^{\natural}_0$,
then $\widehat{X}^{\zeta^{\natural}_0}$ is stationary and
$(\widehat{u}^{\zeta^{\natural}_0},\widehat{X}^{\zeta^{\natural}_0})$
is the optimal couple $(u^{\natural}_0,X^{\natural}_0)$.

\begin{lemma}\label{lemmafinale}
For all $x\in L^2(\Omega)$, $\widehat{u}^x\in\widehat{\U}$
and $\widehat{J}(\widehat{u}^x, x)$ does not depend on $x$.
\end{lemma}
\Dim \noindent Let us consider $X_t^{s,x}$ the solution of equation
\begin{equation*}
\left\lbrace
\begin{array}
[c]{l}
dX^{s,x}_{t}=H_{t}X^{s,x}_tdt+{\displaystyle\sum_{i=1}^{d}}K_{t}^{i}X^{s,x}_{t}dW_{t}^{i}\\
X^{s,x}_s=x,
\end{array}
\right.
\end{equation*}
starting from $x$ at time $s$. We denote, for every $0\leq s\leq
t$, $U(t,s)x:=X_t^{s,x}$.
We notice that
\[
 \widehat{X}^{0,x}_t-\widehat{X}^{0,\zeta^{\natural}_0}_t=x-\zeta^{\natural}_0+
\int_0^sH_s(\widehat{X}^{0,x}_s-\widehat{X}^{0,\zeta^{\natural}_0}_s)ds+\sum_{i=1}^d\int_0^s
K^i_s(\widehat{X}^{0,x}_s-\widehat{X}^{0,\zeta^{\natural}_0}_s)dW^i_s=U(t,0)(x-\zeta^{\natural}_0).
\]
So by the Datko theorem, see e.g. \cite{GM} and \cite{GTinf}, there exist constants $a,C>0$ such that
\[
 \E\vert \widehat{X}^{x}_t-\widehat{X}^{\zeta^{\natural}_0}_t \vert^2
\leq Ce^{-at}\E\vert x-\zeta ^{\natural}_0\vert^2.
\]
So
\[
\E\vert \widehat{X}^{x}\vert^2 \leq Ce^{-at}\E\vert x-\zeta
^{\natural}_0\vert^2+\E\vert
\widehat{X}^{\zeta^{\natural}_0}\vert^2\leq C,
\]
where in the last passage we use that $\widehat{X}^{\zeta^{\natural}}=X^{\natural}$ and it is stationary.

\noindent Again by applying the Datko theorem we obtain
\[
 \lim_{\alpha\rightarrow 0}\alpha \E\int_0^{\infty}e^{-2\alpha s}
 (2\langle SX^{\natural}_s,U(s,0)(x-\zeta^{\natural}_0)\rangle +
 \vert \sqrt{S}U(s,0)(x-\zeta^{\natural}_0)) \vert^2)ds=0.
\]
Moreover
\[
\widehat {u}_t=u^{\natural}_t- \Lambda(t,\overline{P}_t,
\overline{Q}_t)U(0,t)(x-\zeta^{\natural}_0)
\]
It is clear that $u^{\natural} $ belongs to the space of admissible control space  $\widehat{\U}$.

The term
$\tilde{u}_t=\Lambda(t,\overline{P}_t,\overline{Q}_t)U(0,t)(x-\zeta^{\natural}_0)$,
$t \in (0+\infty)$ can be proved to be the optimal control for the
infinite horizon problem with $f=0$ and random initial data
$x-\zeta^{\natural}_0$ :
\[ \inf_{u \in L^2_\P((0,+\infty);\R^k) }\E \int_0^{+\infty} (|\sqrt{S}_sX^u_s|^2+ |u(s)|^2)\,ds. \]
Hence Theorem 5.2 of \cite{GM} can be extended without any
difficulty to get that:
\begin{multline*}
J(0,x-\zeta^{\natural}_0,\tilde{u})  =\E\langle\overline{P}_{0}(x-\zeta^{\natural}_0),x-\zeta^{\natural}_0\rangle+2\E\langle
r_0 ,x-\zeta^{\natural}_0\rangle\\  -\E\int_0 ^{\infty} \vert (
I+{\displaystyle\sum_{i=1}^{d}} \left( D_{t}^{i}\right)
^{\ast}\overline{P}_{t}D_{t}^{i}) ^{-1}(
B_{t}^{*}r_{t}+{\displaystyle\sum_{i=1}^{d}}\left(
D_{t}^{i}\right) ^{*}g_{t}^{i})\vert ^2 ds.
\end{multline*}
Therefore
\[
\E\int_0^{\infty}e^{-2\alpha s}
 |\tilde{u}(s)|^2\,ds \leq  \E\int_0^{\infty}
 |\tilde{u}(s)|^2\,ds \leq C .
\]
This proves that $\widehat {u}$ is an admissible control since it
follows that
\[
 \lim_{\alpha\rightarrow 0}2\alpha \E\int_0^{\infty}e^{-2\alpha s}(\vert u^{\natural}_s\vert^2 -\vert \widehat{u}_s^x \vert^2)ds=0.
\]
\finedim

We can now conclude as follows:
\begin{theorem}\label{teofinale}
 For all $x\in L^2(\Omega)$ the couple $(\widehat{X}^x, \widehat{u}^x)$ is optimal that is
\[
 \widehat{J}(\widehat{u}^x,x)=\min\lbrace\widehat{J}(u,x):u\in\widehat{\U}\rbrace
\]
Moreover the optimal cost, that does not depend on the initial state $x$, is equal to the optimal cost for the periodic (respectively stationary) problem, i.e.
\[
 \widehat{J}(\widehat{u}^x,x)=\overline{J}^{\natural}.
\]
\end{theorem}
\Dim
We denote $ \inf_{u \in \mathcal{U}^{\alpha}} J_\alpha(x,u):=J^*_\alpha(x,u)$.
If $u\in \widehat{\U}$, then for every $\alpha >0$, $u\in\U^{\alpha}$. Consequently for every $\alpha>0$
\[
 2\alpha J_{\alpha}(u,x)\geq 2\alpha J_{\alpha}^{*}.
\]
By taking the limit on both sides we get
\[
 \widehat{J}(x,u)=\underline{\lim}_{\alpha \rightarrow 0}2\alpha J_{\alpha}\geq \underline{\lim}_{\alpha \rightarrow 0}2\alpha J_{\alpha}^{*}=\overline{J}^{\natural}.
\]
By the previous lemma $\widehat{J}(x,\widehat{u}^x)$ is independent on $x$
so we let $x=\zeta^{\natural}_0$, which implies that $\widehat{u}^x=u^{\natural}$ and $\widehat{X}^x=X^{\natural}$. Then
\begin{align*}
& \widehat{J}(\zeta^{\natural}_0,u^{\natural})
=\lim_{\alpha\rightarrow 0}2\alpha \int_0^{\infty}e^{-2\alpha t}
[\vert \sqrt{S}_tX^{\natural}_t\vert^2+\vert u^{\natural}_t\vert^2] dt\\ \nonumber
& = \lim_{\alpha\rightarrow 0}2\alpha (\sum_{k=1}^{+\infty}e^{-2k\alpha})
\int_0^1e^{-2\alpha t}[\vert \sqrt{S}_tX^{\natural}_t\vert^2+\vert u^{\natural}_t\vert^2] dt
=\overline{J}^{\natural},\\ \nonumber
\end{align*}
and this concludes the proof.
\finedim


\begin{thebibliography}{99}


\bibitem{Ben1}
A.~Bensoussan and J.~Frehse.
\newblock {On Bellman Equations of Ergodic Control in $\mathbb{R}^b$}.
\newblock {{\em J. reine angew. Math.} 429(1992),125-160.}


\bibitem{Bi76}
J.-M.~Bismut.
\newblock {Linear quadratic optimal stochastic
control with random coefficients}.
\newblock {\em SIAM J. Contr. Optim.} 14
(1976), 419--444.

\bibitem{Bi78}
J.-M.~Bismut.
\newblock {Contr\^{o}le des syst\`{e}mes lin\'{e}aires quadratiques: applications de l'int\'{e}grale
stochastique.}
\newblock {In {\em S\'{e}minaire de Probabilit\'{e}s, XII}}.
Lecture Notes in Math., 649, Springer, Berlin, 1978.


%
%







%

\bibitem{GM}
G.~Guatteri and F. Masiero.
\newblock{Infinite Horizon and Ergodic Optimal Quadratic Control for an Affine Equation with Stochastic Coefficients,}
\newblock{arxiv:math/0707.0606,{\em Submitted}}

\bibitem{GTinf}
G.~Guatteri and G. Tessitore.
\newblock{Backward Stochastic Riccati Equations and
Infinite Horizon L-Q Optimal Control Problems with Stochastic Coefficients}
\newblock{{\em Applied Mathematics and Optimization,}} 57 (2008), no.2, pp.207-235



\bibitem{KoTa01}
M.~Kohlmann and S.~Tang.
\newblock {  New developments in backward stochastic Riccati
 equations and their applications}.
\newblock {In {\em  Mathematical finance (Konstanz, 2000)}},
Trends Math., Birkh�ser, Basel, 2001.



\bibitem{KoTa02}
M.~Kohlmann and S.~Tang.
\newblock { Global adapted solution of
one-dimensional backward stochastic Riccati equations, with
application to the mean-variance hedging}.
\newblock {\em  Stochastic Process. Appl.} 97
(2002),  1255--288.


\bibitem{KoTa03}
M.~Kohlmann and S.~Tang.
\newblock { Multidimensional backward stochastic Riccati equations and
applications}.
\newblock {\em  SIAM J. Contr. Optim.} 41 (2003),  1696--1721.




%



\bibitem{Peng}
S.~Peng.
\newblock { Stochastic Hamilton-Jacobi-Bellman  Equations.}
\newblock {\em SIAM J. Contr. Optim.} 30 (1992),  284--304.



\bibitem{Tang}
S.~Tang.
\newblock { General linear quadratic optimal control problems with random coefficients:
linear stochastic Hamilton systems and backward stochastic Riccati
equations.}
\newblock {SIAM J. Control Optim. 42 (2003), no. 1, 53--75}




\bibitem{Tess2}
G.Tessitore.
\newblock { Infinite Horizon, Ergodic and Periodic Control for a Stochastic Infinite Dimensional Affine Equation.}
\newblock {\em Journal of Mathematical Systems, Estimation, and Control. 8 n.4 (1998), 1-28}



\end{thebibliography}
\end{document}